\documentclass[11pt, reqno]{amsart}
\usepackage{amsmath, amsthm, amscd, amsfonts, amssymb, color}

\newtheorem{theorem}{Theorem}[section]

\newtheorem{remark}[theorem]{Remark}

\newtheorem{proposition}[theorem]{Proposition}
\newtheorem{definition}{Definition}[section]
\numberwithin{equation}{section}

\def\be{\begin{equation}}
\def\ee{\end{equation}}
\def\ba{\begin{eqnarray}}
\def\ea{\end{eqnarray}}
\def\tilde{\widetilde}

\def\e1{\epsilon}
\def\AAl{\mathcal{A}_{\lambda}}
\def\A0{\stackrel{\circ}{\AAl}}

\def\o1{\omega}
\def\01{\Omega}
\def\c1{\gamma}

\def\g1{\Sigma}
\def\bigcap{\cap}

\def\l1{\Lambda}

\def\v1{\varphi}

\def\d1{\delta}
\def\part{\partial}

\def\f1{\frac}
\def\t1{\theta}
\def\s1{\sqrt{\e1}}
\def\b1{\beta}

\def\bs{\begin{eqnarray*}}
\def\es{\end{eqnarray*}}
\def\m1{\Theta}
\def\w1{\wedge}

\begin{document}
\def\essinf{\mathop{\rm essinf}}

\title[Weighted Morrey estimates for Hausdorff operator and its commutator ]{%
Weighted Morrey estimates for Hausdorff operator and its commutator on the
Heisenberg group}
\author{Jianmiao Ruan, Dashan Fan and Qingyan Wu}
\thanks{The research was supported by the National Natural Science
Foundation of China (Nos. 11471288,\ 11671185,\ 11701250,\ 11771358), the
Zhejiang Provincial Natural Science Foundation of China (No. LY18A010015).}
\address{Jianmiao Ruan, Department of Mathematics\\
Zhejiang International Studies University\\
Hangzhou 310012, China.}
\email{rjmath@163.com}
\address{Dashan Fan, Department of Mathematical Sciences\\
University of Wisconsin-Milwaukee\\
Milwaukee WI 53201, USA.}
\email{fan@uwm.edu}
\address{Qingyan Wu, Department of Mathematics\\
Linyi University\\
Linyi 276005, China.}
\email{wuqingyan@lyu.edu.cn}

\begin{abstract}
In this paper, we study the high-dimensional Hausdorff operators, defined
via a general linear mapping $A$, and their commutators on the weighted
Morrey spaces in the setting of the Heisenberg group. Particularly, under
some assumption on the mapping $A$, we establish their sharp boundedness on
the power weighted Morrey spaces.
\end{abstract}

\maketitle

\section{Introduction.}

\bigskip Let $\mathbf{R}^{n}$ be the Euclidean space of dimension $n.$
Lerner and Liflyand in \cite{LL} studied the Hausdorff operator $H_{\Phi ,A}$
defined by
\begin{equation*}
\text{\ \ }H_{\Phi ,A}(f)(x)=\int_{\mathbf{R}^{n}}\frac{\Phi (y)}{\left\vert
y\right\vert ^{n}}f(xA(y))dy,
\end{equation*}%
where $A(y)$ is an $n\times n$ matrix satisfying $\det A(y)\neq 0$ almost
everywhere in the support of a fixed integrable function $\Phi .$ By
choosing
\begin{equation*}
A(y)=\mathrm{diag}[1/|y|,1/|y|,\dots ,1/|y|],
\end{equation*}%
one then defines $H_{\Phi ,A}$ in this special case by
\begin{equation*}
H_{\Phi }(f)(x)=\int_{\mathbf{R}^{n}}\frac{\Phi (y)}{\left\vert y\right\vert
^{n}}f\left( \frac{x}{\left\vert y\right\vert }\right) dy.
\end{equation*}%
In the definition of $H_{\Phi ,A}(f),$ for simplicity, one may always assume
that functions $f$ \ initially lie in the Schwartz space $S.$ After we
establish the boundedness of $H_{\Phi ,A}(f)$ for $f\in S$ on a normed (or
quasi-normed) space $X,$ we can use a standard dense argument together with
the Hahn-Banach theorem to easily extend the boundedness of \ $H_{\Phi ,A}$
to the whole space $X.$ For the Lebesgue space $L^{p}$ $\left( p\geq
1\right) $ and the Hardy space $H^{1},$ the boundedness of $H_{\Phi }$ (even
$H_{\Phi ,A}$ ) are well established (see \cite{CFL, CZ, Li-2, Li-1, LM-1,
LM-3, Mo, RF, We}). Besides spaces $L^{p}$ and $H^{1}$, the boundedness of $%
H_{\Phi }$ on other function spaces was recently also studied by many
authors (see, for instance, \cite{CFLR, Ka, LM, Mi, RF2, RF3, Wu, WF} and the
references therein). Here, we recommend two recent survey papers \cite{CFW}
and \cite{Li-3} for understanding further the background and historical
development of this research topic. Particularly, it is notable that many
well known operators in analysis can be derived from the Hausdorff operator
if one chooses suitable generating functions $\Phi $ \cite{CFW}$.$

This paper is aimed to study the Hausdorff operator on Morrey spaces. The
classical Morrey spaces introduced by Morrey \cite{Mo1} are a useful work
frame in the study of the existence and regularity of partial differential
equations. It has been obtained that many properties of solutions to partial
differential equations are concerned with the boundedness of some operators
on Morrey type spaces. Therefore, in recent years there has been an
explosion of interest on the boundedness of operators in Morrey type spaces.
For this information, one can refer to \cite{AX, AM,DXY} and references
therein. On the other hand, Chiarenza and Frasca \cite{CF} established the
boundedness of the Hardy-Littlewood maximal operator, of the fractional
integral operator, and of the singular integral operator on Morrey spaces.
Subsequently, Komori and Shirai \cite{KS} extended the results of \cite{CF}
to the weighted Morrey spaces. Alvarez,  Lakey and  Guz\'{m}an-Partida
\cite{ALG} studied the central Morrey spaces.

Inspired by above mentioned research, the purpose of this paper is to study
the boundedness of Hausdorff operator,  as well as its commutator, on the
weighted central Morrey spaces in the setting of the Heisenberg group $\dot{L%
}^{p,\lambda }({\mathbb{H}}^{n};w)$ (see next section for the definition).
We remark that the Hausdorff operator is a linear operator, while we can
view its commutator as a bilinear operator (see Definition 1.1).

The Heisenberg group $\mathbb{H}^{n}$ is a non-commutative nilpotent Lie
group, with the underlying manifold $\mathbf{R}^{2n}\times \mathbf{R}$ and
the group law
\begin{equation*}
x\cdot y=\left( x_{1}+y_{1},x_{2}+y_{2},\cdots
,x_{2n}+y_{2n},x_{2n+1}+y_{2n+1}+2\sum_{j=1}^{n}(y_{j}x_{n+j}-x_{j}y_{n+j})%
\right) ,
\end{equation*}%
where $x=(x_{1},x_{2},\cdots ,x_{2n+1})$, $y=(y_{1},y_{2},\cdots ,y_{2n+1})$%
. The geometric motions on the Heisenberg group $\mathbb{H}^{n}$ are quite
different from those on $\mathbf{R}^{n}$ due to the loss of
interchangeability. On the other hand we find that $\mathbb{H}^{n}$ inherits
some basic structures of $\mathbf{R}^{n}$ . These inheritances are good
enough for us to study the Hausdorff operator on $\mathbb{H}^{n}$. Also,
since  the\ Heisenberg group plays significant roles in many math branches
such as representation theory, several complex analysis, harmonic analysis,
partial deferential equations and quantum mechanics (see \cite{Ho, Sc} for
more details),  an extension of Hausdorff operator to the Heisenberg group
seems interesting and encouraging.

By the definition, the identity element on $\mathbb{H}^{n}$ is $0\in {%
\mathbf{R}}^{2n+1}$, while the inverse element of $x$ is $-x$. The
corresponding Lie algebra is generated by the left-invariant vector fields:
\begin{equation*}
\begin{split}
X_{j}& =\frac{\partial }{\partial x_{j}}+2x_{n+j}\frac{\partial }{\partial
x_{2n+1}},\quad j=1,\cdots ,n, \\
X_{n+j}& =\frac{\partial }{\partial x_{n+j}}-2x_{j}\frac{\partial }{\partial
x_{2n+1}},\quad j=1,\cdots ,n, \\
X_{2n+1}& =\frac{\partial }{\partial x_{2n+1}}.
\end{split}%
\end{equation*}%
The only non-trivial commutator relations are
\begin{equation*}
\lbrack X_{j},X_{n+j}]=-4X_{2n+1},\quad j=1,\cdots ,n.
\end{equation*}

$\mathbb{H}^{n}$ is a homogeneous group in the sense of Folland and Stein
\cite{FS} with dilations
\begin{equation*}
\delta _{r}(x_{1},x_{2},\cdots ,x_{2n},x_{2n+1})=(rx_{1},rx_{2},\cdots
,rx_{2n},r^{2}x_{2n+1}),\quad r>0.
\end{equation*}%
The Haar measure on $\mathbb{H}^{n}$ coincides with the usual Lebesgue
measure on $\mathbf{R}^{2n}\times \mathbf{R}$. We denote the measure of any
measurable set $E\subset \mathbb{H}^{n}$ by $|E|$. It is easy to check that
\begin{equation*}
|\delta _{r}(E)|=r^{Q}|E|,\quad d(\delta _{r}x)=r^{Q}dx.
\end{equation*}%
In the above, $Q=2n+2$ is the homogeneous dimension of $\mathbb{H}^{n}$.

The Heisenberg distance
\begin{equation*}
d(p,q)=d(q^{-1}p,0)=|q^{-1}p|_{h}
\end{equation*}%
is derived from the norm
\begin{equation*}
|x|_{h}=\left[ \left( \sum_{i=1}^{2n}x_{i}^{2}\right) ^{2}+x_{2n+1}^{2}%
\right] ^{\frac{1}{4}},
\end{equation*}%
where $x=(x_{1},x_{2},\cdots ,x_{2n},x_{2n+1}).$

This distance $d$ is left-invariant in the sense that $d(p,q)$ remains
unchanged when $p$ and $q$ are both left-translated by some fixed vector on $%
\mathbb{H}^{n}$. Furthermore, $d$ satisfies the triangular inequality (p.
320 in \cite{KR})
\begin{equation*}
d(p,q)\leq d(p,x)+d(x,q),\quad p,x,q\in \mathbb{H}^{n}.
\end{equation*}

For $r>0$ and $x\in \mathbb{H}^{n}$, the ball and sphere with center $x$ and
radius $r$ on $\mathbb{H}^{n}$ are given by
\begin{equation*}
B(x,r)=\{y\in \mathbb{H}^{n}:d(x,y)<r\},
\end{equation*}%
and
\begin{equation*}
S(x,r)=\{y\in \mathbb{H}^{n}:d(x,y)=r\},
\end{equation*}%
respectively. We know that
\begin{equation*}
|B(x,r)|=|B(0,r)|=\Omega _{Q}r^{Q},
\end{equation*}%
where
\begin{equation}
\Omega _{Q}=\frac{2\pi ^{n+\frac{1}{2}}\Gamma (\frac{n}{2})}{(n+1)\Gamma
(n)\Gamma (\frac{n+1}{2})},  \label{om}
\end{equation}%
is the volume of the unit ball $B(0,1)$ on $\mathbb{H}^{n}$. The area of $%
S(0,1)$ on $\mathbb{H}^{n}$ is $\omega _{Q}=Q\Omega _{Q}$. For more details
about the Heisenberg group one can refer to \cite{FS}.

Now we provide the definition of Hausdorff operators and their commutators
on the Heisenberg group in the following.

\begin{definition}
Let $\Phi $ be a locally integrable function on $\mathbb{H}^{n}$. The Hausdorff
operators on $\mathbb{H}^n$ are defined by
\begin{eqnarray*}
\mathcal{H}_{\Phi}f(x)&=&\int_{\mathbb{H}^n}\frac{\Phi(y)}{|y|_h^{Q}}f\big(%
\delta_{|y|_h^{-1}}x\big)dy, \\
\mathcal{H}_{\Phi ,A}f(x)&=&\int_{\mathbb{H}^{n}}\frac{\Phi (y)}{|y|_{h}^{Q}}%
f\left( A(y)x\right) dy,
\end{eqnarray*}
where $A(y)$ is a matrix-valued function and  $\det A(y)\neq 0$ almost
everywhere in the support of $\Phi $.

If $b\in L_{loc}(\mathbb{H}^n)$, The commutator of Hausdorff operator is
defined by
\begin{equation*}
\mathcal{H}_{\Phi, A}^{b}f=b\mathcal{H}_{\Phi, A}f-\mathcal{H}_{\Phi, A}(bf).
\end{equation*}
\end{definition}

In the above definition, we note that $\mathcal{H}_{\Phi ,A}=\mathcal{H}%
_{\Phi }$ if we choose a special matrix $A$. For a matrix $M$, we will use
the norm $\Vert M\Vert =\sup_{x\in \mathbb{H}^{n},~x\neq 0}|Mx|_{h}/|x|_{h}.$
By Lemma 3.1 in \cite{RFW},
\begin{equation}  \label{introduction-matrix3}
\|M\|^{-Q}\leq |\det M^{-1}|\leq \|M^{-1}\|^Q,
\end{equation}
where $M$ is any invertible $(2n+1)\times (2n+1)$ matrix. Define
\begin{equation*}
\left(\mathcal{H}_{\Phi, A}^{b,1}f\right)(x)=\int_{\|A(y)\|\leq1}\frac{%
\Phi(y)}{|y|_h^{Q}}f\left(A(y)x\right)\left[b(x)-b\left(A(y)x\right)\right]%
dy,
\end{equation*}
\begin{equation*}
\left(\mathcal{H}_{\Phi, A}^{b,2}f\right)(x)=\int_{\|A(y)\|>1}\frac{\Phi(y)}{%
|y|_h^{Q}}f\left(A(y)x\right)\left[b(x)-b\left(A(y)x\right)\right]dy.
\end{equation*}%
It is not difficult to see that the commutator can be rewritten by
\begin{equation*}
\mathcal{H}_{\Phi, A}^{b}f=\mathcal{H}_{\Phi, A}^{b,1}f+\mathcal{H}_{\Phi,
A}^{b,2}f.
\end{equation*}

Here and throughout this paper, we use the notation $A\preceq B$ to denote
that there is a constant $C>0$ independent of all essential values\ and
variables such that $A\leq CB.$ We use the notation $A\simeq B$, if there
exists a positive constant $C$ independent of all essential values and
variables, such that $C^{-1}B\leq A\leq CB$. Also, the class $A_{p}$ denotes
the set of all $A_{p}$ weights whose definition can be found in the next
section.

Now we are in a position to state our results.

\begin{theorem}
\label{th-Aq-weight-Hausdorff}\ \ Let $1\leq p_{1},\ p_{2},\ q<\infty $ and $%
-1/p_{1}\leq \lambda <0$. Suppose that $w\in A_{q}$ with the critical index $%
r_{w}$ for the reverse H\"{o}lder condition. If $p_{1}>p_{2}qr_{w}/(r_{w}-1)$%
, then we have that, for any $1<\delta <r_{w}$,
\begin{equation*}
\Vert {\mathcal{H}}_{\Phi ,A}f\Vert _{\dot{L}^{p_{2},\lambda }({\mathbb{H}}%
^{n};w)}\preceq C_{1}\Vert f\Vert _{\dot{L}^{p_{1},\lambda }({\mathbb{H}}%
^{n};w)},
\end{equation*}%
where
\allowdisplaybreaks{\begin{eqnarray*}
C_{1}&=&\int_{\Vert A(y)\Vert>1}\frac{|\Phi (y)|}{|y|_h^{Q}}\left(\frac{\|A(y)\|^{Q}}{|\det
A (y)|}\right)^{q/p_{1}}\|A(y)\|^{Q\lambda(\delta-1)/\delta}dy \\
&&+\int_{\Vert A(y)\Vert \leq 1}\frac{|\Phi (y)|}{|y|_h^{Q}}\left(\frac{\|A(y)\|^{Q}}{|\det
A (y)|}\right)^{q/p_{1}}\|A(y)\|^{Q\lambda q}dy.
\end{eqnarray*}}
\end{theorem}

\begin{theorem}
\label{th-Aq-weight-commutator}\ \ Let $1\leq p,\ p_{1},\ p_{2},\ q<\infty $
and $-1/p_{1}\leq \lambda <0$. Suppose that $w\in A_{q}$ with the critical
index $r_{w}$ for the reverse H\"{o}lder condition. If $%
1/p>(1/p_{1}+1/p_{2})qr_{w}/(r_{w}-1)$ and $q\leq p_{2}$, then we have that,
for any $1<\delta <r_{w}$,
\begin{equation*}
\Vert {\mathcal{H}}_{\Phi ,A}^{b}f\Vert _{\dot{L}^{p,\lambda }({\mathbb{H}}%
^{n};w)}\preceq C_{2}\Vert f\Vert _{\dot{L}^{p_{1},\lambda }({\mathbb{H}^{n}}%
;w)}\Vert f\Vert _{{CMO}^{p_{2}}({\mathbb{H}}^{n};w)},
\end{equation*}%
where%
\allowdisplaybreaks{\begin{eqnarray*}
C_{2}\!\!\!&=&\!\!\!\int_{\Vert A(y)\Vert>1}\frac{|\Phi (y)|}{|y|_{h}^{Q}}\left(\frac{\|A(y)\|^{Q}}{|\det
A (y)|}\right)^{q/p_{1}}\|A(y)\|^{Q\lambda(\delta-1)/\delta}\max \left\{\frac{\|A(y)\|^{Q}}{|\det
A (y)|},\ \log_{2} \|A(y)\|\right\}dy \\
&&+\int_{\Vert A(y)\Vert \leq 1}\frac{|\Phi (y)|}{|y|_{h}^{Q}}\left(\frac{\|A(y)\|^{Q}}{|\det
A (y)|}\right)^{q/p_{1}}\|A(y)\|^{Q\lambda q}\max \left\{\frac{\|A(y)\|^{Q}}{|\det
A (y)|},\ \log_{2} \frac{1}{\|A(y)\|}\right\}dy.
\end{eqnarray*}}
\end{theorem}

When the weight is reduced to the power function, we have the following enhanced
results.

\begin{theorem}
\label{th-power weight-Hausdorff}\ \ Let $1\leq p<\infty ,\ -1/p\leq \lambda
<0$ and $-Q<\alpha <\infty $. We have that {\
\allowdisplaybreaks{
\begin{eqnarray*}
\Vert {\mathcal{H}}_{\Phi ,A}f\Vert _{\dot{L}^{p, \lambda}(\mathbb{H}^{n};|\cdot|_{h}^{\alpha})}
\preceq C_{3}(\alpha)\Vert f\Vert
_{\dot{L}^{p, \lambda}(\mathbb{H}^{n};|\cdot|_{h}^{\alpha})},
\end{eqnarray*}}} where
\allowdisplaybreaks{\begin{equation*}
C_{3}(\alpha)\!=\!\left\{
\begin{array}{ll}
\!\!\!\displaystyle{\int_{\mathbb{H}^{n}}\frac{|\Phi (y)|}{|y|_{h}^{Q}}\frac{\|A(y)\|^{(Q+\alpha)(\lambda+1/p)}}{\left|\det A(y)\right|^{1/p}}\|A^{-1}(y)\|^{\alpha/p}dy} , & \!\!0<\alpha<\infty, \\
\!\!\!\displaystyle{\int_{\mathbb{H}^{n}}\!\!\frac{|\Phi (y)|}{|y|_{h}^{Q}}\frac{\|A(y)\|^{(Q+\alpha)(\lambda+1/p)}}{\left|\det A(y)\right|^{1/p}}\|A(y)\|^{-\alpha/p}dy}, & \!\!-Q<\alpha\leq0.
\end{array}\right.
\end{equation*}}
\end{theorem}

\begin{theorem}
\label{th-power weight-commutator}\ \ Let $1\leq p,\ p_{1},\ p_{2}<\infty,\
-1/p_{1}\leq\lambda<0$ and $1/p=1/p_{1}+1/p_{2}$.

$\mathrm{(i)}$\ \ If $-Q<\alpha \leq 0$, then we have {\ that
\allowdisplaybreaks{
\begin{eqnarray*}
\Vert {\mathcal{H}}^{b}_{\Phi ,A}f\Vert _{\dot{L}^{p, \lambda}(\mathbb{H}^{n};|\cdot|_{h}^{\alpha})}
\preceq C_{4}\Vert f\Vert
_{\dot{L}^{p_{1}, \lambda}(\mathbb{H}^{n};|\cdot|_{h}^{\alpha})}\Vert b\Vert
_{CMO^{p_{2}}(\mathbb{H}^{n};|\cdot|_{h}^{\alpha})},
\end{eqnarray*}}}where
\allowdisplaybreaks{\begin{equation*}
C_{4}=\int_{\mathbb{H}^{n}}\!\!\frac{|\Phi (y)|}{|y|_{h}^{Q}}\frac{\|A(y)\|^{(Q+\alpha)(\lambda+1/p_{1})}}{\left|\det A(y)\right|^{1/p_{1}}}\|A(y)\|^{-\alpha/p_{1}}\!\max\left\{\!\!\frac{\|A(y)\|^{Q}}{|\det A(y)|}, \left|\log_{2}\|A(y)\|\right|\!\!\right\}\!dy.
\end{equation*}}

$\mathrm{(ii)}$\ \ If $0<\alpha <\infty $ and $p_{2}>(Q+\alpha )/Q$, then we
have {\ that
\allowdisplaybreaks{
\begin{eqnarray*}
\Vert {\mathcal{H}}^{b}_{\Phi ,A}f\Vert _{\dot{L}^{p, \lambda}(\mathbb{H}^{n};|\cdot|_{h}^{\alpha})}
\preceq C_{5}\Vert f\Vert
_{\dot{L}^{p_{1}, \lambda}(\mathbb{H}^{n};|\cdot|_{h}^{\alpha})}\Vert b\Vert
_{CMO^{p_{2}}(\mathbb{H}^{n};|\cdot|_{h}^{\alpha})},
\end{eqnarray*}}}where
\allowdisplaybreaks{\begin{eqnarray*}
C_{5}&=&\int_{\mathbb{H}^{n}}\!\!\frac{|\Phi (y)|}{|y|_{h}^{Q}}\left(\frac{\|A(y)\|^{(Q+\alpha)(\lambda+1/p)}}{\left|\det A(y)\right|^{1/p}}\|A^{-1}(y)\|^{\alpha/p}\right.\\
&&\left.+\frac{\|A(y)\|^{(Q+\alpha)(\lambda+1/p_{1})}}{\left|\det A(y)\right|^{1/p_{1}}}\|A^{-1}(y)\|^{\alpha/p_{1}} \max\left\{\!\!\frac{\|A(y)\|^{Q}}{|\det A(y)|}, \left|\log_{2}\|A(y)\|\right|\!\right\}\right)\!dy.
\end{eqnarray*}}
\end{theorem}

Especially, if $\Vert A^{-1}(y)\Vert $ and $\Vert A(y)\Vert ^{-1}$ are
comparable, the following sharp results hold

\begin{theorem}
\label{th-power weight-Hausdorff-iff}\ \ Let $1\leq p<\infty,\ -1/p\leq
\lambda<0, -Q<\alpha<\infty$ and $\Phi $ be a nonnegative function. Suppose
that there is a constant $C_{0}$ independent of $y$ such that $\Vert
A^{-1}(y)\Vert \leq C_{0}\Vert A(y)\Vert ^{-1}$ for all $y\in
\mathop{\rm
supp}(\Phi )$ . Then ${\mathcal{H}}_{\Phi ,A}$ is bounded on $\dot{L}^{p,\
\lambda}(\mathbb{H}^{n};|\cdot|_{h}^{\alpha})$ if and only if
\begin{equation}  \label{th-sharp1}
\int_{\mathbb{H}^{n}}\frac{\Phi (y)}{|y|_{h}^{Q}}\Vert A(y)\Vert
^{(Q+\alpha)\lambda}dy<\infty .
\end{equation}
\end{theorem}

\begin{theorem}
\label{th-power weight-commutator-iff}\ \ Let $1\leq p,\ p_{1},\
p_{2}<\infty ,\ 1/p=1/p_{1}+1/p_{2},\ -1/p_{1}<\lambda <0$, $-Q<\alpha
<\infty $ and $p_{2}>(Q+\alpha )/Q$ if $0<\alpha <\infty $ or $p_{2}\geq 1$
if $-Q<\alpha \leq 0$. Suppose that $\Phi $ is a nonnegative function and
there is a constant $C_{0}$ independent of $y$ such that $\Vert
A^{-1}(y)\Vert \leq C_{0}\Vert A(y)\Vert ^{-1}$ for all $y\in
\mathop{\rm
supp}(\Phi )$. If $b\in CMO^{p_{2}}(\mathbb{H}^{n};|\cdot |_{h}^{\alpha })$
and \eqref{th-sharp1} holds, then we have\ the following conclusions.

$\mathrm{(i)}$ ${\mathcal{H}}_{\Phi ,A}^{b,1}$ is bounded from $\dot{L}%
^{p_{1},\ \lambda }(\mathbb{H}^{n};|\cdot |_{h}^{\alpha })$ to $\dot{L}^{p,\
\lambda }(\mathbb{H}^{n};|\cdot |_{h}^{\alpha })$ if and only if
\allowdisplaybreaks{\begin{equation*}
\int_{\|A(y)\|\leq 1}\frac{\Phi (y)}{|y|_{h}^{Q}}\Vert A(y)\Vert
^{(Q+\alpha)\lambda}\left|\log_{2}\|A(y)\|\right|dy<\infty.
\end{equation*}}

$\mathrm{(ii)}$ ${\mathcal{H}}_{\Phi ,A}^{b,2}$ is bounded from $\dot{L}%
^{p_{1},\ \lambda }(\mathbb{H}^{n};|\cdot |_{h}^{\alpha })$ to $\dot{L}^{p,\
\lambda }(\mathbb{H}^{n};|\cdot |_{h}^{\alpha })$ if and only if
\allowdisplaybreaks{\begin{equation*}
\int_{\|A(y)\|>1}\frac{\Phi (y)}{|y|_{h}^{Q}}\Vert A(y)\Vert
^{(Q+\alpha)\lambda}\log_{2}\|A(y)\|dy<\infty.
\end{equation*}}
\end{theorem}

Finally in this section, we want to make a few remarks about our main
theorems.

\begin{remark}
\ \ \textrm{Suppose $A(y)=\mathrm{diag}[1/\lambda _{1}(y),\ldots ,1/\lambda
_{2n}(y),1/\lambda_{2n+1}(y)]$ with $\lambda _{i}(y)\neq 0,$ for $i=1,\ldots
,2n+1.$ Denote
\begin{equation*}
\begin{split}
M(y)&=\max\{|\lambda _{1}(y)|,\ldots ,|\lambda _{2n}(y)|,|\lambda
_{2n+1}(y)|^{1/2}\}, \\
m(y)&=\min\{|\lambda _{1}(y)|,\ldots ,|\lambda _{2n}(y)|,|\lambda
_{2n+1}(y)|^{1/2}\}.
\end{split}%
\end{equation*}%
If there is a constant $C\geq 1$ independent of $y$ such that $M(y)\leq
Cm(y) $, then it is easy to check that $A(y)$ satisfies the assumptions of
Theorem \ref{th-power weight-Hausdorff-iff} and Theorem \ref{th-power
weight-commutator-iff}. }
\end{remark}

\begin{remark}
By checking the proof of necessity of Theorem \ref{th-power
weight-commutator-iff}, we find that the necessary condition in (ii) with $%
C_{0}=1$ and in (i) are also true without the assumption \eqref{th-sharp1}.
Therefore, comparing with Theorem \ref{th-power weight-commutator} and
Theorem \ref{th-power weight-commutator-iff}, we raise the following two
questions.

$\mathrm{(i)}$ Do the statements also hold for $\lambda =-1/p_{1}$ in
Theorem \ref{th-power weight-commutator-iff}?

$\mathrm{(ii)}$ Is ${\mathcal{H}}_{\Phi ,A}^{b}$ bounded from $\dot{L}^{p_{1},\
\lambda }(\mathbb{H}^{n};|\cdot |_{h}^{\alpha })$ to $\dot{L}^{p,\ \lambda }(%
\mathbb{H}^{n};|\cdot |_{h}^{\alpha })$ if and only if
\allowdisplaybreaks{\begin{equation*}
\int_{\mathbb{H}^{n}}\frac{\Phi (y)}{|y|_{h}^{Q}}\Vert A(y)\Vert
^{(Q+\alpha)\lambda}\left(1+\left|\log_{2}\|A(y)\|\right|\right)dy<\infty ?
\end{equation*}}
\end{remark}

In the second section, we will introduce some necessary notation and
definitions, as well as some known results to be used later in the paper. We
will prove the main theorems in Section 3.

\section{ Notation and Definitions}

We first recall some standard definitions and notation. The theory of $A_{p}$
weight was first introduced by Muckenhoupt in the Euclidean spaces for
studying the weighted $L^{p}$ boundedness of Hardy-Littlewood maximal
functions in \cite{Mu}. For $A_{p}$ weights on the Heisenberg group one can
refer to \cite{Gu,HPR}. A weight is a nonnegative, locally integrable
function on ${\mathbb{H}}^{n}$.

\begin{definition}
\ Let $1<p<\infty $. We say that a weight $w\in A_{p}(\mathbb{H}^n)$ if
there exists a constant $C$ such that for all balls $B$,
\begin{equation*}
\left( \frac{1}{|B|}\int_{B}w(x)dx\right) \left( \frac{1}{|B|}%
\int_{B}w\left( x\right) ^{-1/(p-1)}dx\right) ^{p-1}\leq C.
\end{equation*}%
We say that a weight $w\in A_{1}(\mathbb{H}^n)$ if there is a constant $C$
such that for all balls $B$,
\begin{equation*}
\frac{1}{|B|}\int_{B}w\left( x\right) dx\leq C\mathop{\rm essinf}%
\limits_{x\in B}w\left( x\right) .
\end{equation*}%
We define
\begin{equation*}
A_{\infty }(\mathbb{H}^n)= \mathop{\cup}\limits_{1\leq p<\infty }A_{p}(%
\mathbb{H}^n).
\end{equation*}
\end{definition}

Following proofs of Propositions 1.4.1, 1.4.2 in \cite{LD} together with the
reverse H\"{o}lder inequality on the Heisenberg group \cite{HPR}, we have
the following results.

\begin{proposition}
\label{th:prop1} $\mathrm{(i)}$ $A_p(\mathbb{H}^n)\subsetneq A_q(\mathbb{H}^n)$%
, for $1\leq p<q<\infty$.\newline
$\mathrm{(ii)}$ If $w\in A_p(\mathbb{H}^n)$, $1<p<\infty$, then there is an $%
\varepsilon>0$ such that $p-\varepsilon>1$ and $w\in A_{p-\varepsilon}(%
\mathbb{H}^n)$.
\end{proposition}

A close relation to $A_{\infty }(\mathbb{H}^n)$ is the reverse H\"{o}lder
condition. If there exist $r>1$ and a fixed constant $C$ such that
\begin{equation*}
\left( \frac{1}{|B|}\int_{B}w(x)^{r}dx\right) ^{1/r}\leq \frac{C}{|B|}%
\int_{B}w(x)dx
\end{equation*}%
for all balls $B\subset \mathbb{H}^{n}$, we then say that $w$ satisfies the
\textit{reverse H\"{o}lder condition of order $r$} and write $w\in RH_{r}(%
\mathbb{H}^n)$. According to Theorem 19 and Corollary 21 in \cite{IMS}, $%
w\in A_{\infty }(\mathbb{H}^n)$ if and only if there exists some $r>1$ such
that $w\in RH_{r}(\mathbb{H}^n)$. Moreover, if $w\in RH_{r}(\mathbb{H}^n),\
r>1$, then $w\in RH_{r+\epsilon }(\mathbb{H}^n)$ for some $\epsilon >0$. We
thus write $r_{w}\equiv \text{sup}\{r>1:w\in RH_{r}(\mathbb{H}^n)\}$ to
denote the \textit{critical index of $w$ for the reverse H\"{o}lder condition%
}.

An important example of $A_{p}(\mathbb{H}^{n})$ weight is the power function
$|x|_{h}^{\alpha }$. By the similar proofs of Propositions 1.4.3 and 1.4.4
in \cite{LD}, we obtain the following properties of power weights.

\begin{proposition}
\label{th:prop2}\ Let $x\in\mathbb{H}^n$. Then\newline
$\mathrm{(i)}$ $|x|_h^{\alpha}\in A_1(\mathbb{H}^n)$ if and only if $%
-Q<\alpha\leq 0$;\newline
$\mathrm{(ii)}$ $|x|_h^{\alpha}\in A_p(\mathbb{H}^n)$, $1<p<\infty$, if and
only if $-Q<\alpha< Q(p-1)$.
\end{proposition}

We will denote by $q_{w}$ the \emph{critical index} for $w$, that is, the
infimum of all the $q$ such that $w$ satisfies the condition $A_{q}$. From
Proposition \ref{th:prop1}, we see that unless $q_{w}=1$, $w$ is never an $%
A_{q_{w}}$ weight. Also by Proposition \ref{th:prop1} and Proposition \ref%
{th:prop2} we see that if $0<\alpha <\infty $, then
\begin{equation}
|x|_{h}^{\alpha }\in \mathop{\bigcap}\limits_{\frac{Q+\alpha }{Q}<p<\infty
}A_{p},  \label{xap}
\end{equation}%
where $(Q+\alpha )/Q$ is the critical index of $|x|_{h}^{\alpha }$.

For any $w\in A_{\infty }(\mathbb{H}^n)$ and any Lebesgue measurable set $E$%
, write $w(E)=\int_{E}w(x)dx$. We have the following standard
characterization of $A_{p}$ weights (see \cite{RFW}).

{%
\allowdisplaybreaks{\begin{proposition}
\label{proAp1}\ \ Let $w\in A_{p}(\mathbb{H}^n)\cap RH_{r}(\mathbb{H}^n), p\geq1$ and $r>1$. Then there
exist constants $C_{1}, C_{2}>0$ such that
\begin{equation*}
C_{1}\left(\frac{|E|}{|B|}\right)^{p}\leq\frac{w(E)}{w\left(B\right)}\leq
C_{2}\left(\frac{|E|}{|B|}\right)^{(r-1)/r}
\end{equation*}
for any measurable subset $E$ of a ball $B$. Especially, for any $\lambda >1$,
\begin{equation*}
w\left( B\left( x_{0},\lambda R\right) \right) \leq C\lambda ^{np}w\left(
B\left( x_{0},R\right) \right) .
\end{equation*}
\end{proposition}}}

\begin{proposition}
\label{proAp2}\ \ If $w\in A_p(\mathbb{H}^n)$, $1\leq p<\infty$, then for
any $f\in L_{loc}^{1}(\mathbb{H}^n)$ and any ball $B\subset \mathbb{H}^n$,
\begin{equation*}
\frac{1}{|B|}\int_B|f(x)|dx\leq C\left(\frac{1}{w(B)}\int_{B}|f(x)|^pw(x)dx%
\right)^{1/p}.
\end{equation*}
\end{proposition}

Given a weight function $w$ on $\mathbb{H}^{n}$, for any measurable set $%
E\subset \mathbb{H}^{n}$, as usual we denote by $L^{p}\left(E; w\right)$ the
weighted Lebesgue space of all functions satisfying
\begin{equation*}
\Vert f\Vert _{L^{p}\left(E; w\right)}=\left(
\int_{E}|f(x)|^{p}w(x)dx\right) ^{1/p}<\infty .
\end{equation*}%
We denote $L^{\infty }(\mathbb{H}^{n};w)=L^{\infty }(\mathbb{H}^{n})$ and $%
\Vert f\Vert _{L^{\infty }(\mathbb{H}^{n};w)}=\Vert f\Vert _{L^{\infty }(%
\mathbb{H}^{n})} $ for $p=\infty $.

\begin{definition}
\label{def-central-Morrey}\ \ Let $1\leq p<\infty ,\ -1/p\leq \lambda <0$
and $w$ be a weight on $\mathbb{H}^{n}$. A function $f\in L_{loc}^{p}(%
\mathbb{H}^{n};w)$ is said to belong to the weighted central Morrey spaces $%
\dot{L}^{p,\lambda }(\mathbb{H}^{n};w)$ if
\begin{equation*}
\Vert f\Vert _{\dot{L}^{p,\lambda }(\mathbb{H}^{n};w)}=\sup_{r>0}\left(
\frac{1}{w(B(0,r))^{1+p\lambda }}\int_{B(0,r)}|f(x)|^{p}w(x)dx\right)
^{1/p}<\infty .
\end{equation*}
\end{definition}

When $\lambda=-1/p$, then $\dot{L}^{p, \lambda}(\mathbb{H}^{n};w)=L^{p}(%
\mathbb{H}^{n};w)$. If $w\equiv1$, one can easily check that $\dot{L}^{p,
\lambda}(\mathbb{H}^{n})$ reduces to ${0}$ when $\lambda<-1/p$.

\begin{definition}
\label{def-CMO}\ \ Let $1\leq p<\infty $ and $w$ be a weight on $\mathbb{H}%
^{n}$. A function $f\in L_{loc}^{p}(\mathbb{H}^{n};w)$ is said to be in the
weighted central BMO spaces $CMO^{p}(\mathbb{H}^{n};w)$ if
\begin{equation*}
\Vert f\Vert _{CMO^{p}(\mathbb{H}^{n};w)}=\sup_{r>0}\left( \frac{1}{w(B(0,r))%
}\int_{B(0,r)}|f(x)-f_{B}|^{p}w(x)dx\right) ^{1/p}<\infty ,
\end{equation*}%
where $f_{B}=\int_{B(0,r)}f(x)dx/|B(0,r)|$.
\end{definition}

The spaces $CMO^{p}(\mathbb{H}^{n};w)$ are quasi-Banach spaces. When $%
1<p<\infty$, then $CMO^{p}(\mathbb{H}^{n};w)$ are Banach spaces after
identifying the functions that differ by a constant almost everywhere. H\"{o}%
lder's inequality shows that $CMO^{p_{1}}(\mathbb{H}^{n};w)\subset
CMO^{p_{2}}(\mathbb{H}^{n};w)$ if $1\leq p_{2}<p_{1}<\infty$. If $w\equiv1$,
we denote the central BMO spaces by $CMO^{p}(\mathbb{H}^{n})$ and we can see
that $BMO(\mathbb{H}^{n})\subset CMO^{p}(\mathbb{H}^{n}),\ 1\leq p<\infty$.

\section{Proof of The Theorems}

In this section, we use $tB(0,r)$ to denote $B(0,tr)$ for any central ball $%
B(0,r)$ in $\mathbb{H}^n$ and any $t>0$. For a $(2n+1)\times (2n+1)$ matrix $%
M$, we denote $MB(0,r)$ the set $\{z\in\mathbb{H}^n\mid z=Mx,~~ x\in
B(0,r)\} $.

\subsection{\noindent \textbf{Proof of Theorem \protect\ref%
{th-Aq-weight-Hausdorff}.}}

By the definition and the Minkowski inequality, {%
\allowdisplaybreaks{\begin{eqnarray}
&&\nonumber \left\|{\mathcal{H}}_{\Phi ,A}f\right\|_{\dot{L}^{p_{2},\lambda}(\mathbb{H}^{n}; w)}\\
\nonumber \!\!&=&\!\!\sup_{r>0}\frac{1}{w(B(0,r))^{\lambda+1/p_{2}}}\left\|\int_{\mathbb{H}^{n}}\frac{\Phi(y)}{|y|_{h}^{Q}}f( A(y)\cdot )dy\right\|_{L^{p_{2}}(B(0,r); w)}\\
\label{proof-weight-Ap-1}\!\!&\leq&\!\!\sup_{r>0}\frac{1}{w(B(0,r))^{\lambda+1/p_{2}}}\int_{\mathbb{H}^{n}}\frac{|\Phi(y)|}{|y|_{h}^{Q}}\|f( A(y)\cdot )\|_{L^{p_{2}}(B(0,r); w)}dy.
\end{eqnarray}}}Since $p_{1}>p_{2}qr_{w}/(r_{w}-1)=p_{2}qr_{w}^{\prime }$,
there is $1<\gamma <r_{w}$ such that $p_{1}/q=p_{2}\gamma ^{\prime
}=p_{2}\gamma /(\gamma -1)$. In view of the H\"{o}lder inequality and the
reverse H\"{o}lder condition, we obtain that {%
\allowdisplaybreaks{\begin{eqnarray*}
&&\|f( A(y)\cdot)\|_{L^{p_{2}}(B(0,r); w)}\\
&\preceq& \left(\int_{B(0,r)}|f(A(y)x)|^{p_{1}/q}dx\right)^{q/p_{1}}\left(\int_{B(0,r)}w(x)^{\gamma}dx\right)^{1/(\gamma p_{2})}\\
&\preceq&|\det A^{-1}(y)|^{q/p_{1}}\left(\int_{A(y)B(0,r)}|f(x)|^{p_{1}/q}dx\right)^{q/p_{1}}\left(\int_{B(0,r)}w(x)^{\gamma}dx\right)^{1/(\gamma p_{2})}\\
&\preceq& |\det A^{-1}(y)|^{q/p_{1}}|B(0,r)|^{-q/p_{1}}w(B(0,r))^{1/p_{2}}\left(\int_{A(y)B(0,r)}|f(x)|^{p_{1}/q}dx\right)^{q/p_{1}}.
\end{eqnarray*}}}Proposition \ref{proAp2} and \eqref{introduction-matrix3}
show that {%
\allowdisplaybreaks{\begin{eqnarray*}
\nonumber & &\left(\int_{A(y)B(0,r)}|f(x)|^{p_{1}/q}dx\right)^{q/p_{1}}\\
&\leq& \left(\int_{B(0,\|A(y)\|r)}|f(x)|^{p_{1}/q}dx\right)^{q/p_{1}}\\
\nonumber&\preceq&\|A(y)\|^{Qq/p_{1}}|B(0,r)|^{q/p_{1}}\left(\frac{1}{w(B(0,\|A(y)\|r))}\int_{B(0,\|A(y)\|r)}|f(x)|^{p_{1}}w(x)dx\right)^{1/p_{1}},
\end{eqnarray*}}}which implies that {%
\allowdisplaybreaks{\begin{eqnarray}
& &\|f(A(y)\cdot )\|_{L^{p_{2}}((B(0,r); w)}\label{proof-weight-Ap-2}\\
\nonumber &\preceq&|\det A^{-1}(y)|^{q/p_{1}}\|A(y)\|^{Qq/p_{1}}\frac{w(B(0,r))^{1/p_{2}}}{w(B(0,\|A(y)\|r)^{1/p_{1}}}\|f\|_{L^{p_{1}}\left(B(0,\|A(y)\|r);w\right)}.
\end{eqnarray}}}Therefore, we infer from \eqref{proof-weight-Ap-1} and %
\eqref{proof-weight-Ap-2} that {%
\allowdisplaybreaks{\begin{eqnarray}
\label{proof-weight-Ap-3}&& \|{\mathcal{H}}_{\Phi ,A}f\|_{\dot{L}^{p_{2},\lambda}(\mathbb{H}^{n}; w)}\\
\nonumber&\preceq& \|f\|_{\dot{L}^{p_{1},\lambda}(\mathbb{H}^{n}; w)} \sup_{r>0}\int_{\mathbb{H}^{n}}\frac{|\Phi(y)|}{|y|_{h}^{Q}}\left(\frac{\|A(y)\|^{Q}}{|\det A(y)|}\right)^{q/p_{1}}\left(\frac{w(B(0,\|A(y)\|r)}{w(B(0,r))}\right)^{\lambda}dy.
\end{eqnarray}}}

If $\|A(y)\|>1$, Proposition \ref{proAp1} shows that, for any $%
1<\delta<r_{w} $, {%
\allowdisplaybreaks{\begin{eqnarray}
\label{proof-weight-Ap-4}\frac{w(B(0,\|A(y)\|r)}{w(B(0,r))}\succeq \left(\frac{|B(0,\|A(y)\|r)|}{|B(0,r)|}\right)^{(\delta-1)/\delta}=\|A(y)\|^{Q(\delta-1)/\delta}.
\end{eqnarray}}} If $\|A(y)\|\leq1$, by Proposition \ref{proAp1} again, we
have {%
\allowdisplaybreaks{\begin{eqnarray}
\label{proof-weight-Ap-5}\frac{w(B(0,\|A(y)\|r)}{w(B(0,r))}\succeq \left(\frac{|B(0,\|A(y)\|r)|}{|B(0,r)|}\right)^{q}=\|A(y)\|^{Qq}.
\end{eqnarray}}}Thus we complete the proof of Theorem \ref%
{th-Aq-weight-Hausdorff} by \eqref{proof-weight-Ap-3}-\eqref{proof-weight-Ap-5}. $\hfill \Box $%
\newline

\subsection{\noindent \textbf{Proof of Theorem \protect\ref%
{th-Aq-weight-commutator}.}}

By the definition,
\allowdisplaybreaks{\begin{eqnarray}
\nonumber &&\|{\mathcal{H}}^{b}_{\Phi, A}f\|_{\dot{L}^{p, \lambda}(\mathbb{H}^n; w)}\\
\nonumber&\preceq& \sup_{r>0}\frac{1}{w(B(0,r))^{\lambda+1/p}}\left(\int_{B(0,r)}\left|\int_{\mathbb{H}^n}\frac{\Phi(y)}{|y|_{h}^{Q}}f(A(y)x)\left[b(x)-b_{B(0,r)}\right]dy\right|^{p}w(x)dx\right)^{1/p}\\
\nonumber&&+\sup_{r>0}\frac{1}{w(B(0,r))^{\lambda+1/p}}\left(\int_{B(0,r)}\left|\int_{\mathbb{H}^n}\frac{\Phi(y)}{|y|_{h}^{Q}}f(A(y)x)\left[b_{B(0,r)}-b_{A(y)B(0,r)}\right]dy\right|^{p}w(x)dx\right)^{1/p}\\
\nonumber&&+\sup_{r>0}\frac{1}{w(B(0,r))^{\lambda+1/p}}\left(\int_{B(0,r)}\left|\int_{\mathbb{H}^n}\frac{\Phi(y)}{|y|_{h}^{Q}}f(A(y)x)\left[b_{A(y)B(0,r)}-b(A(y)x)\right]dy\right|^{p}w(x)dx\right)^{1/p}\\
\label{proof-th-Aq-commutator-1}&:=& I+II+III.
\end{eqnarray}} H\"{o}lder's inequality and Theorem \ref%
{th-Aq-weight-Hausdorff} show that
\allowdisplaybreaks{\begin{eqnarray}
\nonumber I&=&\sup_{r>0}\frac{1}{w(B(0,r))^{\lambda+1/p}}\left(\int_{B(0,r)}\left|{\mathcal{H}}_{\Phi, A}f(x)\right|^{p}|b(x)-b_{B(0,r)}|^{p}w(x)dx\right)^{1/p}\\
\nonumber &\preceq&\|{\mathcal{H}}_{\Phi, A}f\|_{\dot{L}^{p_{3}, \lambda}(\mathbb{H}^n; w)}\|b\|_{CMO^{p_{2}}(\mathbb{H}^n; w)}\\
\nonumber &\preceq& \|f\|_{\dot{L}^{p_{1}, \lambda}(\mathbb{H}^n; w)}\|b\|_{CMO^{p_{2}}(\mathbb{H}^n; w)}\left(\int_{\Vert A(y)\Vert>1}\frac{|\Phi (y)|}{|y|_{h}^{Q}}\left(\frac{\|A(y)\|^{Q}}{|\det
A (y)|}\right)^{q/p_{1}}\|A(y)\|^{Q\lambda(\delta-1)/\delta}dy \right.\\
\label{proof-th-Aq-commutator-2}&&+\left.\int_{\Vert A(y)\Vert \leq 1}\frac{|\Phi (y)|}{|y|_{h}^{Q}}\left(\frac{\|A(y)\|^{Q}}{|\det
A (y)|}\right)^{q/p_{1}}\|A(y)\|^{Q\lambda q}dy\right),
\end{eqnarray}}where $1/p=1/p_{2}+1/p_{3}$.

According to the Minkowski inequality,
\allowdisplaybreaks{\begin{eqnarray*}
II\!\!&\preceq&\!\!\sup_{r>0}\frac{1}{w(B(0,r))^{\lambda+1/p}}\int_{\mathbb{H}^n}\frac{|\Phi(y)|}{|y|_{h}^{Q}}\left|b_{B(0,r)}-b_{A(y)B(0,r)}\right|\left(\int_{B(0,r)}\left|f(A(y)x)\right|^{p}w(x)dx\right)^{1/p}\!dy.
\end{eqnarray*}}By the similar argument as Theorem \ref%
{th-Aq-weight-Hausdorff}, we have, for any $1<\delta<r_{w}$,
\allowdisplaybreaks{\begin{eqnarray}
\nonumber&&\sup_{r>0}\frac{1}{w(B(0,r))^{\lambda+1/p}}\|f(A(y)\cdot)\|_{L^{p}(B(0,r);w)}\\
\label{proof-th-Aq-commutator-3}&\preceq &\|f\|_{\dot{L}^{p_{1},\lambda}(B(0,r);w)}\left\{\begin{array}{ll}
\!\!\!\displaystyle{\left(\frac{\|A(y)\|^{Q}}{|\det A(y)|}\right)^{q/p_{1}}\|A(y)\|^{Q\lambda(\delta-1)/\delta}} ,& \!\!\|A(y)\|>1, \\
\!\!\!\displaystyle{\left(\frac{\|A(y)\|^{Q}}{|\det A(y)|}\right)^{q/p_{1}}\|A(y)\|^{Q\lambda q}}, &\!\! \|A(y)\|\leq 1.
\end{array}\right.
\end{eqnarray}}

Next, we estimate $\sup_{r>0}\left\vert b_{B(0,r)}-b_{A(y)B(0,r)}\right\vert
$. If $\Vert A(y)\Vert >1$, there exists a nonnegtive integer $k_{0}\geq 0$
satisfying
\begin{equation*}
2^{k_{0}}<\Vert A(y)\Vert \leq 2^{k_{0}+1}.
\end{equation*}%
Therefore
\allowdisplaybreaks{\begin{eqnarray}
\nonumber&&\left|b_{B(0,r)}-b_{A(y)B(0,r)}\right|\\
\nonumber&\preceq&\left|b_{B(0,r)}-b_{2B(0,r)}\right|+\left|b_{2B(0,r)}-b_{4B(0,r)}\right|+\dots+\left|b_{2^{k_{0}}B(0,r)}-b_{2^{k_{0}+1}B(0,r)}\right|\\
\label{proof-th-Aq-commutator-4}&&+\left|b_{2^{k_{0}+1}B(0,r)}-b_{A(y)B(0,r)}\right|.
\end{eqnarray}}Since for any $k\in \mathbf{Z}$, Proposition \ref{proAp2}
implies that
\allowdisplaybreaks{\begin{eqnarray}
\nonumber &&\left|b_{2^{k}B(0,r)}-b_{2^{k+1}B(0,r)}\right|\\
&\leq& \frac{2^{Q}}{|2^{k+1}B(0,r)|}\int_{2^{k+1}B(0,r)}\left|b(x)-b_{2^{k+1}B(0,r)}\right|dx\\
\nonumber &\preceq&\left(\frac{1}{w(2^{k+1}B(0,r))}\int_{2^{k+1}B(0,r)}\left|b(x)-b_{2^{k+1}B(0,r)}\right|^{q}w(x)dx\right)^{1/q}\\
\nonumber &\preceq&\left(\frac{1}{w(2^{k+1}B(0,r))}\int_{2^{k+1}B(0,r)}\left|b(x)-b_{2^{k+1}B(0,r)}\right|^{p_{2}}w(x)dx\right)^{1/_{p_{2}}}\\
\label{proof-th-Aq-commutator-5}&\preceq& \|b\|_{CMO^{p_{2}}(\mathbb{H}^n;w)},
\end{eqnarray}}where the third inequality is achieved by H\"{o}lder's
inequality and $q\leq p_{2}$.

On the other hand,
\allowdisplaybreaks{\begin{eqnarray}
\nonumber &&\left|b_{2^{k_{0}+1}B(0,r)}-b_{A(y)B(0,r)}\right|\\
\nonumber&\preceq&\frac{1}{|A(y)B(0,r)|}\int_{A(y)B(0,r)}\left|b(x)-b_{2^{k_{0}+1}B(0,r)}\right|dx\\
\nonumber&\preceq& \frac{2^{Q(k_{0}+1)}}{|\det A(y)|}\frac{1}{|2^{k_{0}+1}B(0,r)|}\int_{2^{k_{0}+1}B(0,r)}\left|b(x)-b_{2^{k_{0}+1}B(0,r)}\right|dx\\
\nonumber&\preceq& \frac{\|A(y)\|^{Q}}{|\det A(y)|}\left(\frac{1}{w(2^{k_{0}+1}B(0,r))}\int_{2^{k_{0}+1}B(0,r)}\left|b(x)-b_{2^{k_{0}+1}B(0,r)}\right|^{q}w(x)dx\right)^{1/q}\\
\label{proof-th-Aq-commutator-6}&\preceq&\frac{\|A(y)\|^{Q}}{|\det A(y)|}\|b\|_{CMO^{p_{2}}(\mathbb{H}^n;w)}.
\end{eqnarray}}Then, it follows from \eqref{proof-th-Aq-commutator-4}-%
\eqref{proof-th-Aq-commutator-6} that, for $\|A(y)\|>1$,
\allowdisplaybreaks{\begin{eqnarray}
\nonumber \left|b_{B(0,r)}-b_{A(y)B(0,r)}\right|&\preceq& \left(k_{0}+1+\frac{\|A(y)\|^{Q}}{|\det A(y)|}\right)\|b\|_{CMO^{p_{2}}(\mathbb{H}^n;w)}\\
\label{proof-th-Aq-commutator-7}&\preceq& \max\left\{\log_{2}\|A(y)\|, \frac{\|A(y)\|^{Q}}{|\det A(y)|}\right\}\|b\|_{CMO^{p_{2}}(\mathbb{H}^n;w)}.
\end{eqnarray}}Similar to the proceeding argument, for $\|A(y)\|\leq1$,
\allowdisplaybreaks{\begin{eqnarray}
\label{proof-th-Aq-commutator-8} \left|b_{B(0,r)}-b_{A(y)B(0,r)}\right|
 \preceq \max\left\{\log_{2}\frac{1}{\|A(y)\|}, \frac{\|A(y)\|^{Q}}{|\det A(y)|}\right\}\|b\|_{CMO^{p_{2}}(\mathbb{H}^n;w)}.
\end{eqnarray}}Therefore \eqref{proof-th-Aq-commutator-3}, %
\eqref{proof-th-Aq-commutator-7} and \eqref{proof-th-Aq-commutator-8} yield
that
\allowdisplaybreaks{\begin{eqnarray}
\label{proof-th-Aq-commutator-9} II\!\!&\preceq&\!\!\|f\|_{\dot{L}^{p_{1}, \lambda}(\mathbb{H}^n; w)}\|b\|_{CMO^{p_{2}}(\mathbb{H}^n; w)}\\
\nonumber &&\!\!\!\times\!\!\left(\!\int_{\Vert A(y)\Vert>1}\!\!\frac{|\Phi (y)|}{|y|_{h}^{Q}}\left(\!\frac{\|A(y)\|^{Q}}{|\det
A (y)|}\!\right)^{q/p_{1}}\!\!\!\!\|A(y)\|^{Q\lambda(\delta-1)/\delta}\!\max\left\{\!\log_{2}\|A(y)\|, \frac{\|A(y)\|^{Q}}{|\det A(y)|}\!\right\}\!dy \right.\\
\nonumber&&\!\!\!+\!\!\left.\int_{\Vert A(y)\Vert \leq 1}\!\!\frac{|\Phi (y)|}{|y|_{h}^{Q}}\left(\!\frac{\|A(y)\|^{Q}}{|\det
A (y)|}\!\right)^{q/p_{1}}\!\!\!\!\|A(y)\|^{Q\lambda q}\!\max\left\{\!\log_{2}\frac{1}{\|A(y)\|}, \frac{\|A(y)\|^{Q}}{|\det A(y)|}\!\right\}\!dy\!\!\right).
\end{eqnarray}} Now we turn to estimate the term $III$. Using the Minkowski
inequality again,
\allowdisplaybreaks{\begin{eqnarray*}
III\!\!&\preceq&\!\!\sup_{r>0}\frac{1}{w(B(0,r))^{\lambda+1/p}}\int_{\mathbb{H}^n}\frac{|\Phi(y)|}{|y|_{h}^{Q}}\left\|f(A(y)\cdot)\left[b_{A(y)B(0,r)}-b(A(y)\cdot)\right]\right\|_{L^{p}(B(0,r);w)}dy.
\end{eqnarray*}}On the other hand, since $1/p>(1/p_{1}+1/p_{2})qr^{%
\prime}_{w}$, we can choose $p_{4},\ p_{5}$ satisfying $1/p_{4}>qr^{%
\prime}_{w}/p_{1}, 1/p_{5}>qr^{\prime}_{w}/p_{2}$ and $1/p=1/p_{4}+1/p_{5}$.
Then \eqref{proof-th-Aq-commutator-3} implies that,
\allowdisplaybreaks{\begin{eqnarray}
 \nonumber&&\sup_{r>0}\frac{1}{w(B(0,r))^{\lambda+1/p}}\|f(A(y)\cdot)[b_{A(y)B(0,r)}-b(A(y)\cdot)]\|_{L^{p}(B(0,r);w)}\\
\nonumber&\preceq& \sup_{r>0}\frac{1}{w(B(0,r))^{\lambda+1/p}}\|f(A(y)\cdot)\|_{L^{p_{4}}(B(0,r);w)}\|b_{A(y)B(0,r)}-b(A(y)\cdot)\|_{L^{p_{5}}(B(0,r);w)}\\
\nonumber&=&\sup_{r>0}\frac{1}{w(B(0,r))^{\lambda+1/p_{4}}}\|f(A(y)\cdot)\|_{L^{p_{4}}(B(0,r);w)}\\
\nonumber&&\times\left(\frac{1}{w(B(0,r))}\int_{B(0,r)}\left|b_{A(y)B(0,r)}-b(A(y)x)\right|^{p_{5}}w(x)dx\right)^{1/p_{5}}\\
\nonumber&\preceq&\sup_{r>0}\left(\frac{1}{w(B(0,r))}\int_{B(0,r)}\left|b_{A(y)B(0,r)}-b(A(y)x)\right|^{p_{5}}w(x)dx\right)^{1/p_{5}}\\
\label{proof-th-Aq-commutator-10}&&\times\|f\|_{\dot{L}^{p_{1},\lambda}(\mathbb{H}^n;w)}\left\{\begin{array}{ll}
\!\!\!\displaystyle{\left(\frac{\|A(y)\|^{Q}}{|\det A(y)|}\right)^{q/p_{1}}\|A(y)\|^{Q\lambda(\delta-1)/\delta}} ,& \!\!\|A(y)\|>1, \\
\!\!\!\displaystyle{\left(\frac{\|A(y)\|^{Q}}{|\det A(y)|}\right)^{q/p_{1}}\|A(y)\|^{Q\lambda q}}, &\!\! \|A(y)\|\leq 1.
\end{array}\right.
\end{eqnarray}}It follows from Proposition \ref{proAp2} and $q\leq p_{2}$
that
\allowdisplaybreaks{\begin{eqnarray}
\nonumber&&\left(\frac{1}{w(B(0,r))}\int_{B(0,r)}\left|b_{A(y)B(0,r)}-b(A(y)x)\right|^{p_{5}}w(x)dx\right)^{1/p_{5}}\\
\nonumber&\preceq&\frac{1}{w(B(0,r))^{1/p_{5}}}\left(\int_{B(0,r)}\left|b(A(y)x)-b_{B(0,\|A(y)\|r)}\right|^{p_{5}}w(x)dx\right)^{1/p_{5}}\\
\nonumber&&+\left|b_{B(0,\|A(y)\|r)}-b_{A(y)B(0,r)}\right|\\
\nonumber&\preceq&\frac{1}{w(B(0,r))^{1/p_{5}}}\left(\int_{B(0,r)}\left|b(A(y)x)-b_{B(0,\|A(y)\|r)}\right|^{p_{5}}w(x)dx\right)^{1/p_{5}}\\
\nonumber&&+
\frac{\|A(y)\|^{Q}}{|\det A(y)|}\left(\frac{1}{w(B(0,\|A(y)\|r))}\int_{B(0,\|A(y)\|r)}\left|b(x)-b_{B(0,\|A(y)\|r)}\right|^{p_{2}}w(x)dx\right)^{1/p_{2}}\\
\nonumber&\preceq& \frac{1}{w(B(0,r))^{1/p_{5}}}\left(\int_{B(0,r)}\left|b(A(y)x)-b_{B(0,\|A(y)\|r)}\right|^{p_{5}}w(x)dx\right)^{1/p_{5}}\\
\nonumber&&+\frac{\|A(y)\|^{Q}}{|\det A(y)|}\|b\|_{CMO^{p_{2}}(\mathbb{H}^n;w)}.
\end{eqnarray}}Since $1/p_{5}>qr^{\prime}_{w}/p_{2}$, there exists $%
1<\beta<r_{w}$ such that $1/p_{5}=q\beta^{\prime}/p_{2}=q\beta/\left(p_{2}(%
\beta-1)\right)$. Due to the H\"{o}lder inequality and the reverse H\"{o}%
lder condition, we have
\allowdisplaybreaks{\begin{eqnarray}
\nonumber &&\frac{1}{w(B(0,r))^{1/p_{5}}}\left(\int_{B(0,r)}\left|b(A(y)x)-b_{B(0,\|A(y)\|r)}\right|^{p_{5}}w(x)dx\right)^{1/p_{5}}\\
\nonumber &\preceq& \frac{1}{w(B(0,r))^{1/p_{5}}}\left(\int_{B(0,r)}\left|b(A(y)x)-b_{B(0,\|A(y)\|r)}\right|^{p_{2}/q}dx\right)^{q/p_{2}}\left(\int_{B(0,r)}w(x)^{\beta}dx\right)^{1/(\beta p_{5})}\\
\nonumber &\preceq& \frac{1}{w(B(0,r))^{1/p_{5}}\left|\det A(y)\right|^{q/p_{2}}} \left(\int_{B(0,\|A(y)\|r)}\left|b(x)-b_{B(0,\|A(y)\|r)}\right|^{p_{2}/q}dx\right)^{q/p_{2}}\\
\nonumber&& \times \left(\int_{B(0,r)}w(x)^{\beta}dx\right)^{1/(\beta p_{5})}\\
\nonumber &\preceq& \frac{1}{w(B(0,r))^{1/p_{5}}}|B(0,r)|^{(1-\beta)/(\beta p_{5})}\left(\int_{B(0,r)}w(x)dx\right)^{1/p_{5}}\left(\frac{|B(0,\|A(y)\|r)|}{|\det A(y)|}\right)^{q/p_{2}}\\
\nonumber&&\times\left(\frac{1}{w(B(0,\|A(y)\|r))}\int_{B(0,\|A(y)\|r)}\left|b(x)-b_{B(0,\|A(y)\|r)}\right|^{p_{2}}w(x)dx\right)^{1/p_{2}}\\
\nonumber&\preceq&\left(\frac{\|A(y)\|^{Q}}{|\det A(y)|}\right)^{q/p_{2}}\left(\frac{1}{w\left(B(0,\|A(y)\|r)\right)}\int_{B(0,\|A(y)\|r)}\left|b(x)-b_{B(0,\|A(y)\|r)}\right|^{p_{2}}w(x)dx\right)^{1/p_{2}}\\
\nonumber&\preceq&\left(\frac{\|A(y)\|^{Q}}{|\det A(y)|}\right)^{q/p_{2}} \|b\|_{CMO^{p_{2}}(\mathbb{H}^n;w)}\\
\nonumber &\preceq&\frac{\|A(y)\|^{Q}}{|\det A(y)|} \|b\|_{CMO^{p_{2}}(\mathbb{H}^n;w)},
\end{eqnarray}}which yields that
\allowdisplaybreaks{\begin{eqnarray}
\nonumber&&\sup_{r>0}\left(\frac{1}{w(B(0,r))}\int_{B(0,r)}\left|b_{A(y)B(0,r)}-b(A(y)x)\right|^{p_{5}}w(x)dx\right)^{1/p_{5}}\\
\label{proof-th-Aq-commutator-11}&\preceq& \frac{\|A(y)\|^{Q}}{|\det A(y)|} \|b\|_{CMO^{p_{2}}(\mathbb{H}^n;w)}.
\end{eqnarray}}Then, we infer from \eqref{proof-th-Aq-commutator-10} and %
\eqref{proof-th-Aq-commutator-11} that
\allowdisplaybreaks{\begin{eqnarray}
\nonumber III\!\!&\preceq&\|f\|_{\dot{L}^{p_{1}, \lambda}(\mathbb{H}^n; w)}\|b\|_{CMO^{P_{2}}(\mathbb{H}^n; w)}\\
\nonumber &&\times\left(\int_{\Vert A(y)\Vert>1}\frac{|\Phi (y)|}{|y|_{h}^{Q}}\left(\frac{\|A(y)\|^{Q}}{|\det
A (y)|}\right)^{1+q/p_{1}}\|A(y)\|^{Q\lambda(\delta-1)/\delta}dy \right.\\
\label{proof-th-Aq-commutator-12}&&+\left.\int_{\Vert A(y)\Vert \leq 1}\frac{|\Phi (y)|}{|y|_{h}^{Q}}\left(\frac{\|A(y)\|^{Q}}{|\det
A (y)|}\right)^{1+q/p_{1}}\|A(y)\|^{Q\lambda q}dy\right).
\end{eqnarray}}By combining \eqref{proof-th-Aq-commutator-1}, %
\eqref{proof-th-Aq-commutator-2} , \eqref{proof-th-Aq-commutator-9} and %
\eqref{proof-th-Aq-commutator-12}, we finish the proof of the theorem. $%
\hfill\Box$\newline

Here and in after, we some time use $w(\cdot )$ for $|\cdot |_{h}^{\alpha }$
for the sake of convenience.

\subsection{\noindent \textbf{Proof of Theorem \protect\ref{th-power
weight-Hausdorff}. }}

Similar to the proof of the Theorem \ref{th-Aq-weight-Hausdorff}, we have {%
\allowdisplaybreaks{\begin{eqnarray}
&&\nonumber \left\|{\mathcal{H}}_{\Phi ,A}f\right\|_{\dot{L}^{p,\lambda}(\mathbb{H}^n; w)}\\
\label{proof-power weight-1}\!\!&\preceq&\!\!\sup_{r>0}\frac{1}{w(B(0,r))^{\lambda+1/p}}\int_{\mathbb{H}^n}\frac{|\Phi(y)|}{|y|_{h}^{Q}}\|f( A(y)\cdot )\|_{L^{p}(B(0,r); w)}dy.
\end{eqnarray}}}
A changing of variables $A(y)x=z$ yields that
\begin{eqnarray}
&& \left\Vert f(A(y)x)\right\Vert _{L^{p}(B(0,r);|x|^{\alpha })}
\label{proof-power weight-2} \\
&\preceq &\left\{
\begin{array}{ll}
|\det A^{-1}(y)|^{1/p}\Vert A^{-1}(y)\Vert ^{\alpha /p}\left(
\int_{B(0,\|A(y)\|r)}\left\vert f(z)\right\vert ^{p}|z|_h^{\alpha}dz\right)
^{1/p}, & \!\!0<\alpha<\infty, \\
|\det A^{-1}(y)|^{1/p}\Vert A(y)\Vert ^{-\alpha /p}\left(
\int_{B(0,\|A(y)\|r)}\left\vert f(z)\right\vert^{p}|z|_h^{\alpha}dz\right)
^{1/p}, & \!\!-Q<\alpha\leq0. \\
\end{array}%
\right.  \notag
\end{eqnarray}
For simplicity, we just consider the case $0<\alpha<\infty$, since the proof
of the case $-Q<\alpha\leq0$ is essentially similar. Combining %
\eqref{proof-power weight-1},\eqref{proof-power weight-2} and the simple
fact
\begin{eqnarray}  \label{proof-power weight-3}
w(B(0,r))=\int_{|x|_h< r}|x|_h^{\alpha}dx=\frac{\omega_{Q}}{Q+\alpha}
r^{Q+\alpha},\ \ -Q<\alpha<\infty,
\end{eqnarray}%
where $\omega_{Q}$ is the surface area of the unit sphere in $\mathbb{H}^n$,
we have {%
\allowdisplaybreaks{\begin{eqnarray*}
&&\left\|{\mathcal{H}}_{\Phi ,A}f\right\|_{\dot{L}^{p,\lambda}(\mathbb{H}^n; w)}\\
\!\!&\preceq&\!\!\|f\|_{\dot{L}^{p,\lambda}(\mathbb{H}^n; w)}\sup_{r>0}\int_{\mathbb{H}^n}\frac{|\Phi(y)|}{|y|_{h}^{Q}}|\det A^{-1}(y)|^{1/p}\Vert A^{-1}(y)\Vert ^{\alpha /p}\left(\frac{w(B(0,\|A(y)\|r))}{w(B(0,r))}\right)^{\lambda+1/p}dy\\
\!\!&\preceq&\!\!\|f\|_{\dot{L}^{p,\lambda}(\mathbb{H}^n; w)}\int_{\mathbb{H}^n}\frac{|\Phi(y)|}{|y|_{h}^{Q}}|\det A^{-1}(y)|^{1/p}\Vert A^{-1}(y)\Vert ^{\alpha /p}\|A(y)\|^{(\lambda+1/p)(Q+\alpha)}dy.
\end{eqnarray*}}}This proves the theorem. $\hfill\Box$\newline

\subsection{\noindent \textbf{Proof of Theorem \protect\ref{th-power
weight-commutator}. }}

According to the same argument as Theorem \ref{th-Aq-weight-commutator}, we
have
\allowdisplaybreaks{\begin{eqnarray}
\nonumber \|{\mathcal{H}}^{b}_{\Phi, A}f\|_{\dot{L}^{p, \lambda}(\mathbb{H}^n; w)}\preceq I+II+III,
\end{eqnarray}}here $I,\ II,\ III$ are the same as in %
\eqref{proof-th-Aq-commutator-1}.

We first consider the case that $0<\alpha<\infty$. The H\"{o}lder inequality
and Theorem \ref{th-power weight-Hausdorff} show that
\allowdisplaybreaks{\begin{eqnarray}
\nonumber I&\preceq &\|{\mathcal{H}}_{\Phi, A}f\|_{\dot{L}^{p_{1}, \lambda}(\mathbb{H}^n; w)}\|b\|_{CMO^{p_{2}}(\mathbb{H}^n; w)}\\
\nonumber &\preceq & \|f\|_{\dot{L}^{p_{1}, \lambda}(\mathbb{H}^n; w)}\|b\|_{CMO^{p_{2}}(\mathbb{H}^n; w)}\\
\label{proof-power weight-commtator-1}&&\times \int_{\mathbb{H}^n}\frac{|\Phi (y)|}{|y|_{h}^{Q}}\frac{\|A(y)\|^{(Q+\alpha)(\lambda+1/p_{1})}}{\left|\det A(y)\right|^{1/p_{1}}}\|A^{-1}(y)\|^{\alpha/p_{1}}dy.
\end{eqnarray}}

By the Minkowski inequality and H\"{o}lder's inequality again,
\allowdisplaybreaks{\begin{eqnarray}
\nonumber II&\preceq &\!\!\sup_{r>0}\frac{1}{w(B(0,r))^{\lambda+1/p}}\int_{\mathbb{H}^n}\frac{|\Phi(y)|}{|y|_{h}^{Q}}\left|b_{B(0,r)}-b_{A(y)B(0,r)}\right|\left(\int_{B(0,r)}\left|f(A(y)x)\right|^{p}|x|_h^{\alpha}dx\right)^{1/p}\!dy\\
\nonumber &\preceq &\!\!\sup_{r>0}\frac{1}{w(B(0,r))^{\lambda+1/p_{1}}}\int_{\mathbb{H}^n}\frac{|\Phi(y)|}{|y|_{h}^{Q}}\left|b_{B(0,r)}-b_{A(y)B(0,r)}\right|\left(\int_{B(0,r)}\left|f(A(y)x)\right|^{p_{1}}|x|_h^{\alpha}dx\right)^{1/p_{1}}\!dy\\
 \label{proof-power weight-commtator-2}&\preceq &\|f\|_{\dot{L}^{p_{1}, \lambda}(\mathbb{H}^n; w)}\\
\nonumber&&\times \sup_{r>0}\int_{\mathbb{H}^n}\frac{|\Phi (y)|}{|y|_{h}^{Q}}\frac{\|A(y)\|^{(Q+\alpha)(\lambda+1/p_{1})}}{\left|\det A(y)\right|^{1/p_{1}}}\|A^{-1}(y)\|^{\alpha/p_{1}}\left|b_{B(0,r)}-b_{A(y)B(0,r)}\right|dy,
\end{eqnarray}}where the last inequality is achieved by a similar argument
as in the proof of Theorem \ref{th-power weight-Hausdorff}. A similar
discussion as in \eqref{proof-th-Aq-commutator-4}-%
\eqref{proof-th-Aq-commutator-8} shows that
\allowdisplaybreaks{\begin{eqnarray}
\label{proof-power weight-commtator-3} \sup_{r>0}\left|b_{B(0,r)}-b_{A(y)B(0,r)}\right| \preceq \max\left\{\left|\log_{2}\|A(y)\| \right|, \frac{\|A(y)\|^{Q}}{|\det A(y)|}\right\}\|b\|_{CMO^{p_{2}}(\mathbb{H}^n;w)},
\end{eqnarray}}where $p_{2}>(Q+\alpha )/Q$. Then we infer from %
\eqref{proof-power weight-commtator-2} and
\eqref{proof-power
weight-commtator-3} that
\allowdisplaybreaks{\begin{eqnarray}
\label{proof-power weight-commtator-3-1}II&\preceq &\|f\|_{\dot{L}^{p_{1}, \lambda}(\mathbb{H}^n; w)}\|b\|_{CMO^{p_{2}}(\mathbb{H}^n;w)}\\
\nonumber &&\times \int_{\mathbb{H}^n}\!\!\!\frac{|\Phi(y)|}{|y|_{h}^{Q}}\frac{\|A(y)\|^{(Q+\alpha)(\lambda+1/p_{1})}}{\left|\det A(y)\right|^{1/p_{1}}}\|A^{-1}(y)\|^{\alpha/p_{1}}\max\left\{\!\!\left|\log_{2}\|A(y)\| \right|\!, \! \frac{\|A(y)\|^{Q}}{|\det A(y)|}\!\!\right\}\!dy.
\end{eqnarray}}

It is not difficult to check that
\allowdisplaybreaks{\begin{eqnarray*}
III\!\!&\preceq&\!\!\sup_{r>0}\frac{1}{w(B(0,r))^{\lambda+1/p}}\int_{\mathbb{H}^n}\frac{|\Phi(y)|}{|y|_{h}^{Q}}\|f(A(y)\cdot)[b_{A(y)B(0,r)}-b(A(y)\cdot)]\|_{L^{p}(B(0,r);w)}dy,
\end{eqnarray*}}and for any $r>0$,
\allowdisplaybreaks{\begin{eqnarray}
 \nonumber&&\frac{1}{w(B(0,r))^{\lambda+1/p}}\|f(A(y)\cdot)[b_{A(y)B(0,r)}-b(A(y)\cdot)]\|_{L^{p}(B(0,r);w)}\\
\nonumber&\preceq&\frac{1}{w(B(0,r))^{\lambda+1/p_{1}}}\|f(A(y)\cdot)\|_{L^{p_{1}}(B(0,r);w)}\\
\nonumber&&\times\left(\frac{1}{w(B(0,r))}\int_{B(0,r)}\left|b_{A(y)B(0,r)}-b(A(y)x)\right|^{p_{2}}|x|_h^{\alpha}dx\right)^{1/p_{2}}\\
\nonumber&\preceq&\|f\|_{\dot{L}^{p_{1}, \lambda}(\mathbb{H}^n; w)}\frac{\|A(y)\|^{(Q+\alpha)(\lambda+1/p_{1})}}{\left|\det A(y)\right|^{1/p_{1}}}\|A^{-1}(y)\|^{\alpha/p_{1}}\\
\label{proof-power weight-commtator-4}&&\times\left(\frac{1}{w(B(0,r))}\int_{B(0,r)}\left|b_{A(y)B(0,r)}-b(A(y)x)\right|^{p_{2}}|x|_h^{\alpha}dx\right)^{1/p_{2}}.
\end{eqnarray}} On the other hand, the Minkowski inequality and Proposition %
\ref{proAp2} imply that
\allowdisplaybreaks{\begin{eqnarray}
\nonumber&&\left(\frac{1}{w(B(0,r))}\int_{B(0,r)}\left|b_{A(y)B(0,r)}-b(A(y)x)\right|^{p_{2}}|x|_h^{\alpha}dx\right)^{1/p_{2}}\\
\nonumber&\leq&\left(\frac{1}{w(B(0,r))}\int_{B(0,r)}\left|b(A(y)x)-b_{B(0,\|A(y)\|r)}\right|^{p_{2}}|x|_h^{\alpha}dx\right)^{1/p_{2}}+\left|b_{A(y)B(0,r)}-b_{B(0,\|A(y)\|r)}\right|\\
\nonumber&\preceq&\|A^{-1}(y)\|^{\alpha/p_{2}}\left|\det A^{-1}(y)\right|^{1/p_{2}}\left(\frac{1}{w(B(0,r))}\int_{A(y)B(0,r)}\left|b(z)-b_{B(0,\|A(y)\|r)}\right|^{p_{2}}|z|^{\alpha}dz \right)^{1/p_{2}}\\
\nonumber&&+\frac{1}{|A(y)B(0,r)|}\int_{A(y)B(0,r)}\left|b(z)-b_{B(0,\|A(y)\|r)}\right|dz\\
\nonumber&\preceq&\|A^{-1}(y)\|^{\alpha/p_{2}}\left|\det A^{-1}(y)\right|^{1/p_{2}}\left(\frac{w(B(0,\|A(y)\|r))}{w(B(0,r))}\right)^{1/p_{2}}\|b\|_{CMO^{p_{2}}(\mathbb{H}^n;w)}\\
\nonumber&&+\frac{\|A(y)\|^{Q}}{|\det A(y)|}\|b\|_{CMO^{p_{2}}(\mathbb{H}^n;w)}\\
\label{proof-power weight-commtator-5}&\preceq&\left(\|A^{-1}(y)\|^{\alpha/p_{2}}\left|\det A^{-1}(y)\right|^{1/p_{2}}\|A(y)\|^{(Q+\alpha)/p_{2}}+\frac{\|A(y)\|^{Q}}{|\det A(y)|}\right)\|b\|_{CMO^{p_{2}}(\mathbb{H}^n;w)}.
\end{eqnarray}} We infer from \eqref{proof-power weight-commtator-4}, %
\eqref{proof-power weight-commtator-5} and \eqref{introduction-matrix3} that
\allowdisplaybreaks{\begin{eqnarray}
\nonumber III\!\!&\preceq&\!\!\|f\|_{\dot{L}^{p_{1}, \lambda}(\mathbb{H}^n; w)}\|b\|_{CMO^{p_{2}}(\mathbb{H}^n;w)}\int_{\mathbb{H}^n}\frac{|\Phi(y)|}{|y|_{h}^{Q}}\left(\frac{\|A^{-1}(y)\|^{\alpha/p}\|A(y)\|^{(Q+\alpha)(\lambda+1/p)}}{|\det A(y)|^{1/p}}\right.\\
\label{proof-power weight-commtator-6}\!\!&&\!\! \left.+\frac{\|A(y)\|^{(Q+\alpha)(\lambda+1/p_{1})}\|A^{-1}(y)\|^{\alpha/p_{1}}}{\left|\det A(y)\right|^{1/p_{1}}}\frac{\|A(y)\|^{Q}}{|\det A(y)|}\right)dy.
\end{eqnarray}} Thus we complete the proof of the case $0<\alpha<\infty$ by %
\eqref{proof-power weight-commtator-1},
\eqref{proof-power
weight-commtator-3-1} and \eqref{proof-power weight-commtator-6}.

Next we consider the case that $-Q<\alpha \leq 0$. By the previously used
argument, we have
\allowdisplaybreaks{\begin{eqnarray}
\nonumber I&\preceq & \|f\|_{\dot{L}^{p_{1}, \lambda}(\mathbb{H}^n; w)}\|b\|_{CMO^{p_{2}}(\mathbb{H}^n; w)}\\
\label{proof-power weight-commtator-7}&&\times \int_{\mathbb{H}^n}\frac{|\Phi (y)|}{|y|_{h}^{Q}}\frac{\|A(y)\|^{(Q+\alpha)(\lambda+1/p_{1})}}{\left|\det A(y)\right|^{1/p_{1}}}\|A(y)\|^{-\alpha/p_{1}}dy\\
\nonumber&=&\|f\|_{\dot{L}^{p_{1}, \lambda}(\mathbb{H}^n; w)}\|b\|_{CMO^{p_{2}}(\mathbb{H}^n; w)}
\int_{\mathbb{H}^n}\frac{|\Phi (y)|}{|y|_{h}^{Q}}\frac{\|A(y)\|^{(Q+\alpha)\lambda+Q/p_{1}}}{\left|\det A(y)\right|^{1/p_{1}}}dy.
\end{eqnarray}}
\allowdisplaybreaks{\begin{eqnarray}
\label{proof-power weight-commtator-8}II&\preceq &\|f\|_{\dot{L}^{p_{1}, \lambda}(\mathbb{H}^n; w)}\|b\|_{CMO^{p_{2}}(\mathbb{H}^n;w)}\\
\nonumber &&\times \int_{\mathbb{H}^n}\!\!\!\frac{|\Phi(y)|}{|y|_{h}^{Q}}\frac{\|A(y)\|^{(Q+\alpha)\lambda+Q/p_{1}}}{\left|\det A(y)\right|^{1/p_{1}}}\max\left\{\!\!\left|\log_{2}\|A(y)\| \right|\!, \! \frac{\|A(y)\|^{Q}}{|\det A(y)|}\!\!\right\}\!dy.
\end{eqnarray}}
\allowdisplaybreaks{\begin{eqnarray}
\nonumber III\!\!&\preceq&\!\!\|f\|_{\dot{L}^{p_{1}, \lambda}(\mathbb{H}^n; w)}\|b\|_{CMO^{p_{2}}(\mathbb{H}^n;w)}\int_{\mathbb{H}^n}\frac{|\Phi(y)|}{|y|_{h}^{Q}}\left(\frac{\|A(y)\|^{(Q+\alpha)\lambda+Q/p}}{|\det A(y)|^{1/p}}\right.\\
\nonumber\!\!&&\!\! \left.+\frac{\|A(y)\|^{(Q+\alpha)\lambda+Q/p_{1}}}{\left|\det A(y)\right|^{1/p_{1}}}\frac{\|A(y)\|^{Q}}{|\det A(y)|}\right)dy\\
\label{proof-power weight-commtator-9}\!\!&\preceq&\!\!\|f\|_{\dot{L}^{p_{1}, \lambda}(\mathbb{H}^n; w)}\|b\|_{CMO^{p_{2}}(\mathbb{H}^n;w)}\int_{\mathbb{H}^n}\frac{\|A(y)\|^{(Q+\alpha)\lambda+Q/p_{1}}}{\left|\det A(y)\right|^{1/p_{1}}}\frac{\|A(y)\|^{Q}}{|\det A(y)|}dy,
\end{eqnarray}}where the last inequality is obtained by %
\eqref{introduction-matrix3}. Thus we finish the proof of the theorem.$%
\hfill \Box $\newline

\subsection{\noindent \textbf{Proof of Theorem \protect\ref{th-power
weight-Hausdorff-iff}.}\ }

\bigskip

If $\Vert A^{-1}(y)\Vert \preceq \Vert A(y)\Vert ^{-1}$, then %
\eqref{introduction-matrix3} gives that
\begin{equation}
\Vert A^{-1}(y)\Vert ^{Q}\simeq \Vert A(y)\Vert ^{-Q}\simeq |\det A^{-1}(y)|.
\label{proof-H(i)-in-9}
\end{equation}%
The \textquotedblleft $if$\textquotedblright\ part of Theorem \ref{th-power
weight-Hausdorff-iff} is easily obtained from Theorem \ref{th-power
weight-Hausdorff}. Next we will show the \textquotedblleft $only\ if$%
\textquotedblright\ part.

We just consider the case that $-1/p<\lambda<0$, since the theorem is
exactly the Corollary 1.6 in \cite{RFW} if $\lambda=-1/p$. Let $%
f^{*}(x)=|x|_h^{(Q+\alpha)\lambda}$. By the fact of
\eqref{proof-power
weight-3}, it is not difficult to check that
\begin{eqnarray*}
\|f^{*}\|_{\dot{L}^{p,\lambda}(\mathbb{H}^n;|\cdot|_{h}^{\alpha})}=\frac{%
(Q+\alpha)^{\lambda}}{\omega_{Q}^{\lambda}(1+p\lambda)^{1/p}}.
\end{eqnarray*}%
Simple calculation shows that, for $\Phi\geq 0$, {%
\allowdisplaybreaks{\begin{eqnarray*}
\left({\mathcal{H}}_{\Phi,A}f^{*}\right)(x)&=&\int_{\mathbb{H}^n}\frac{\Phi(y)}{|y|_{h}^{Q}}|A(y)x|^{(Q+\alpha)\lambda}dy\\
&\geq& f^{*}(x)\int_{\mathbb{H}^n}\frac{\Phi(y)}{|y|_{h}^{Q}}\|A(y)\|^{(Q+\alpha)\lambda}dy,
\end{eqnarray*}}}which completes the proof of Theorem \ref{th-power
weight-Hausdorff-iff}.$\hfill\Box$\newline

\subsection{\noindent \textbf{Proof of Theorem \protect\ref{th-power
weight-commutator-iff}.}\ }

\bigskip

The sufficient part is easily obtained by \eqref{proof-H(i)-in-9} and
Theorem \ref{th-power weight-commutator}. It now remains to prove the
necessary part.

(i) Let $b^{*}(x)=\ln |x|_h$. A simple calculation tells us that
\begin{equation*}
b^{*}_{B(0,r)}=\frac{1}{|B(0,r)|}\int_{B(0,r)}\ln|x|_hdx=\ln r-1/Q,
\end{equation*}
which implies that
\allowdisplaybreaks{\begin{eqnarray*}
\int_{B(0,r)}\left|b^{*}(x)-b^{*}_{B(0,r)}\right|^{p_{2}}|x|_h^{\alpha}dx&=&\omega_{Q}\int_{0}^{r}\left|1/Q-\ln(r/\rho)\right|^{p_{2}}\rho^{\alpha+Q-1}d\rho\\
&=&\omega_{Q}r^{Q+\alpha}\int_{1}^{\infty}\left|1/Q-\ln t\right|^{p_{2}}t^{-(Q+\alpha+1)}dt.
\end{eqnarray*}}Therefore,
\allowdisplaybreaks{\begin{eqnarray*}
\|b^{*}\|_{CMO^{p_{2}}(\mathbb{H}^n;|\cdot|_{h}^{\alpha})}=c (Q+\alpha)^{1/p_{2}},
\end{eqnarray*}}where
\begin{equation*}
c=\left(\int_{1}^{\infty}\left|1/Q-\ln
t\right|^{p_{2}}t^{-(Q+\alpha+1)}dt\right)^{1/p_{2}}<\infty.
\end{equation*}
Let $f^{*}(x)=|x|_h^{(Q+\alpha)\lambda}$. Then $f^{*}\in \dot{L}%
^{p_{1},\lambda}(\mathbb{H}^n;|\cdot|_{h}^{\alpha})\cap \dot{L}^{p,\lambda}(%
\mathbb{H}^n;|\cdot|_{h}^{\alpha})$ and
\begin{eqnarray}  \label{proof-power-weight-commutator-iff-1}
\|f^{*}\|_{\dot{L}^{p_{1},\lambda}(\mathbb{H}^n;|\cdot|_{h}^{\alpha})}=\frac{%
(Q+\alpha)^{\lambda}}{\omega_{Q}^{\lambda}(1+p_{1}\lambda)^{1/p_{1}}},\ \
\|f^{*}\|_{\dot{L}^{p,\lambda}(\mathbb{H}^n;|\cdot|_{h}^{\alpha})}=\frac{%
(Q+\alpha)^{\lambda}}{\omega_{Q}^{\lambda}(1+p\lambda)^{1/p}}.
\end{eqnarray}%
By definition of the commutator of the Hausdorff operator,
\allowdisplaybreaks{\begin{eqnarray*}
\left({\mathcal{H}}^{b^{*},1}_{\Phi,A}f^{*}\right)(x)=\int_{\|A(y)\|\leq1}\frac{\Phi(y)}{|y|_{h}^{Q}}|A(y)x|_h^{(Q+\alpha)\lambda}\ln\left(\frac{|x|_h}{|A(y)x|_h}\right)dy.
\end{eqnarray*}} Since $\|A(y)\|\leq 1$,
\allowdisplaybreaks{\begin{eqnarray*}
\ln\left(\frac{|x|_h}{|A(y)x|_h}\right)\geq \ln\left(\frac{1}{\|A(y)\|}\right)\geq 0.
\end{eqnarray*}}This inequality implies that
\allowdisplaybreaks{\begin{eqnarray*}
\|{\mathcal{H}}^{b^{*},1}_{\Phi,A}f^{*}\|_{\dot{L}^{p,\lambda}(\mathbb{H}^n;|\cdot|_{h}^{\alpha})}\geq \|f^{*}\|_{\dot{L}^{p,\lambda}(\mathbb{H}^n;|\cdot|_{h}^{\alpha})}\int_{\|A(y)\|\leq 1}\frac{\Phi(y)}{|y|_{h}^{Q}}\|A(y)\|^{(Q+\alpha)\lambda}\ln\left(\frac{1}{\|A(y)\|}\right)dy.
\end{eqnarray*}}Therefore, the boundedness of $H^{b^{*},1}_{\Phi,A}$ from $%
\dot{L}^{p_{1},\lambda}(\mathbb{H}^n;|\cdot|_{h}^{\alpha})$ to $\dot{L}%
^{p,\lambda}(\mathbb{H}^n;|\cdot|_{h}^{\alpha})$ and %
\eqref{proof-power-weight-commutator-iff-1} show that
\begin{equation*}
\int_{\|A(y)\|\leq 1}\frac{\Phi(y)}{|y|_{h}^{Q}}\|A(y)\|^{(Q+\alpha)\lambda}%
\ln\left(\frac{1}{\|A(y)\|}\right)dy<\infty.
\end{equation*}

(ii) Let $\widetilde{b}(x)=\ln\left(1/|x|_h\right)$. Then $\widetilde{b}\in
CMO^{p_{2}}(\mathbb{H}^n;|\cdot|_{h}^{\alpha})$ and
\begin{equation*}
\|\widetilde{b}\|_{CMO^{p_{2}}(\mathbb{H}^n;|\cdot|_{h}^{\alpha})}=\|b^{*}%
\|_{CMO^{p_{2}}(\mathbb{H}^n;|\cdot|_{h}^{\alpha})},
\end{equation*}
where $b^{*}$ is as in (i). Without loss of generality, we assume that the
constant $C_{0}>1$, since the case of $C_{0}=1$ is easier to deal with.
Taking $f^{*}$ be as in (i), we have
\allowdisplaybreaks{\begin{eqnarray}\label{proof-th-power weight-commutator-iff}
\left({\mathcal{H}}^{\tilde{b},2}_{\Phi,A}f^{*}\right)(x)
=\int_{\|A(y)\|>1}\frac{\Phi(y)}{|y|_h^{Q}}|A(y)x|_h^{(Q+\alpha)\lambda}\ln\left(\frac{|A(y)x|_h}{|x|_h}\right)dy=I_{1}+I_{2},
\end{eqnarray}}where
\allowdisplaybreaks{\begin{eqnarray*}
I_{1}=\int_{\|A(y)\|> C_{0}}\frac{\Phi(y)}{|y|_h^{Q}}|A(y)x|_h^{(Q+\alpha)\lambda}\ln\left(\frac{|A(y)x|_h}{|x|_h}\right)dy,\\
I_{2}=\int_{1<\|A(y)\|\leq C_{0}}\frac{\Phi(y)}{|y|_h^{Q}}|A(y)x|_h^{(Q+\alpha)\lambda}\ln\left(\frac{|A(y)x|_h}{|x|_h}\right)dy.
\end{eqnarray*}}Since $\|A(y)\|^{-1}\leq\|A^{-1}(y)\|\leq C_{0}\|A(y)\|^{-1}$%
, if $\|A(y)\|> C_{0}$,
\allowdisplaybreaks{\begin{eqnarray}\label{proof-th-power weight-commutator-iff-1}
\ln\left(\frac{|A(y)x|_h}{|x|_h}\right)\geq \ln\left(\frac{1}{\|A^{-1}(y)\|}\right)\geq\ln\left(\frac{\|A(y)\|}{C_{0}}\right)\geq 0,
\end{eqnarray}} and if $1<\|A(y)\|\leq C_{0}$,
\allowdisplaybreaks{\begin{eqnarray}\label{proof-th-power weight-commutator-iff-2}
\left|\ln\left(\frac{|A(y)x|_h}{|x|_h}\right)\right|\leq \ln C_{0}.
\end{eqnarray}}It follows from \eqref{proof-th-power weight-commutator-iff-1}
and \eqref{proof-th-power weight-commutator-iff-2} that
\allowdisplaybreaks{\begin{eqnarray*}
I_{1}\geq |x|_h^{(Q+\alpha)\lambda}\int_{\|A(y)\|\geq C_{0}}\frac{\Phi(y)}{|y|_h^{Q}}\|A(y)\|^{(Q+\alpha)\lambda}\ln\left(\frac{\|A(y)\|}{C_{0}}\right)dy\geq0,\\
|I_{2}|\leq |x|_h^{(Q+\alpha)\lambda}\ln C_{0}\int_{1<\|A(y)\|<C_{0}}\frac{\Phi(y)}{|y|_h^{Q}}\|A(y)\|^{(Q+\alpha)\lambda}dy.
\end{eqnarray*}}These two inequalities and
\eqref{proof-th-power
weight-commutator-iff} tell us that
\allowdisplaybreaks{\begin{eqnarray*}
&&\|f^{*}\|_{\dot{L}^{p,\lambda}(\mathbb{H}^n;|\cdot|_h^{\alpha})}\int_{\|A(y)\|> C_{0}}\frac{\Phi(y)}{|y|_h^{Q}}\|A(y)\|^{(Q+\alpha)\lambda}\ln\left(\frac{\|A(y)\|}{C_{0}}\right)dy\\
&\leq&\|{\mathcal{H}}^{\tilde{b},2}_{\Phi,A}f^{*}\|_{\dot{L}^{p,\lambda}(\mathbb{H}^n;|\cdot|_h^{\alpha})}+\|f^{*}\|_{\dot{L}^{p,\lambda}(\mathbb{H}^n;|\cdot|_h^{\alpha})} \ln C_{0}\int_{1<\|A(y)\|\leq C_{0}}\frac{\Phi(y)}{|y|_h^{Q}}\|A(y)\|^{(Q+\alpha)\lambda}dy.
\end{eqnarray*}}Therefore, the boundedness of $H^{\widetilde{b},2}_{\Phi,A}$
from $\dot{L}^{p_{1},\lambda}(\mathbb{H}^n;|\cdot|_{h}^{\alpha})$ to $\dot{L}%
^{p,\lambda}(\mathbb{H}^n;|\cdot|_{h}^{\alpha})$ and \eqref{th-sharp1}, %
\eqref{proof-power-weight-commutator-iff-1} show that
\begin{equation*}
\int_{\|A(y)\|> C_{0}}\frac{\Phi(y)}{|y|_h^{Q}}\|A(y)\|^{(Q+\alpha)\lambda}%
\ln\left(\frac{\|A(y)\|}{C_{0}}\right)dy<\infty.
\end{equation*}
Using \eqref{th-sharp1} again, we finish the proof of Theorem \ref{th-power
weight-commutator-iff}.$\hfill\Box$\newline


\end{document}